\documentclass[11pt,reqno]{amsart}
\usepackage{amsmath,amssymb,amsthm,upref,graphicx,mathrsfs}

\usepackage{color,xcolor}
\usepackage[
  colorlinks=true,
  linkcolor=blue,
  citecolor=blue,
  urlcolor=blue]{hyperref}

\newtheorem{thm}{Theorem}[section]
\newtheorem{cor}[thm]{Corollary}

\newtheorem{prop}[thm]{Proposition}

\numberwithin{equation}{section}

\begin{document}

\leftline{ \scriptsize}

\vspace{1.3 cm}
\title
{Structured Singular values on some generalized stochastic matrices}
\author{Kijti Rodtes $^{\ast}$ and Mutti Ur Rehman Abbasi $^\circ$}
\thanks{{\scriptsize
		\newline MSC(2010): 15A18, 15B51, 15A60.  \\ Keywords: Structured singular values, Generalized doubly stochastic matrices, Norms of matrices.  \\
			$^{*}$ Department of Mathematics, Faculty of Science, Naresuan University, Phitsanulok 65000, Thailand, (kijtir@nu.ac.th). \\
		$^{\circ}$ Department of Mathematics, Sukkur IBA University, 65200, Sukkur-Pakistan, (mutti.rehman@iba-suk.edu.pk) and Department of Chemical Engineering, MIT, Cambridge, Massachusetts 02139, USA, (muttiur@mit.edu).
		\\}}
\hskip -0.4 true cm

\maketitle


\begin{abstract}In this paper, the exact values of the structured singular values of some generalized stochastic complex matrices is explicit in term of the constant row (column) sum.    
\end{abstract}

\vskip 0.2 true cm


\pagestyle{myheadings}
\markboth{\rightline {\scriptsize Kijti Rodtes and Mutti Ur Rehman Abbasi}}
{\leftline{\scriptsize }}
\bigskip
\bigskip


\vskip 0.4 true cm

\section{Introduction}

The structured singular value is a useful tool for analyzing the robustness of stability and performance of dynamical system, \cite{DoyWS}.  This notion was introduced originally by John Doyle in \cite{Doy}.  In fact, it is a map $\mu_\mathbb{B}: M_n(\mathbb{C})\longrightarrow [0,\infty)$ defined by
$$\mu_\mathbb{B}(M):= \begin{cases} 0, &\mbox{if } \det(I+M\Delta)\neq 0, \forall \Delta\in \mathbb{B}\\
(\min_{\Delta\in \mathbb{B}}\{\|\Delta\|_2\;|\; \det(I+M\Delta)=0\})^{-1},& \mbox{otherwise. }  \end{cases} $$
Here, $M_n(\mathbb{C})$ denotes the set of all $n\times n$ complex matrices, $\mathbb{B}$ denotes a set of all diagonal blocks $$\mathbb{B}:=\{\operatorname{diag}(\delta_1I_{k_1},\dots,\delta_rI_{k_r},\Delta_1,\dots,\Delta_s)\,|\;\delta_i's\in \mathbb{C}, \Delta_j's\in M_{n_i}(\mathbb{C})\}  \subseteq M_n(\mathbb{C})$$
with $\sum_{i=1}^rk_i+\sum_{j=1}^sn_j=n$ and $\|A\|_2$ denotes the matrix $2$-norm of $A$ (which is also equal to the largest singular value of $A$, $\sigma_{\max}(A),$ for all $A\in M_n(\mathbb{C})$).
Note that the complex perturbation block structure $\mathbb{B}$ depends on the values of $r$ (number of repeated scalar blocks), $s$ (number of full blocks),  $k_i$'s and $n_i$'s (dimensions of each block).

Some basic properties of the structured singular value are deduced from the definition  including (see (\cite{Doy}): for any block structure $\mathbb{B}\subseteq M_n(\mathbb{C})$,
 \begin{enumerate}
 	\item $\mu_\mathbb{B}(\alpha M)=|\alpha| \mu_\mathbb{B}(M)$, for any $M\in M_n(\mathbb{C})$ and $\alpha\in \mathbb{C}$. 
 	\item $\mu_\mathbb{B}(M)=\sigma_{\max}(M)$, for any $M\in\mathbb{B}$.
 	\item $\rho(M)\leq \mu_\mathbb{B}(M)\leq\sigma_{\max}(M)$, for any $M\in M_n(\mathbb{C})$, where $\rho(M)$ denotes the spectral radius (=magnitude of the largest eigenvalue) of $M$.
 	\item Let $\mathcal{D}_{\mathbb{B}}:=\{\operatorname{diag}(d_1I_{k_1},\dots,d_rI_{k_r},d_{r+1}I_{n_1},\dots,d_{r+s}I_{n_s})\;|\;d_i\in (0,\infty)\}$ and denote $\mathcal{U}_{\mathbb{B}}:=\{U\in \mathbb{B}\;|\; U \hbox{ is unitary}\}$.  Then,  for any $M\in M_n(\mathbb{C})$,
 	\begin{equation}\label{propineqmu1}
\max_{U\in \mathcal{U}_{\mathbb{B}}}\rho(MU)\leq \mu_\mathbb{B}(M)\leq \inf_{D\in\mathcal{D}_{\mathbb{B}} }\sigma_{\max}(DMD^{-1}).
 	\end{equation}
\end{enumerate} 
The first inequality in (\ref{propineqmu1}) is shown to be equality in \cite{Doy}; namely,
\begin{equation}\label{lowerexact}
\mu_{\mathbb{B}}(M)=\max_{U\in \mathcal{U}_\mathbb{B}}\rho(MU).  
\end{equation} 
Unfortunately, this optimization problem is not convex.  Precisely, the function $l(U):=\rho(MU)$ can have multiple local maxima which are not global and so the direct calculation of $\max_{U\in \mathcal{U}_\mathbb{B}}\rho(MU)$ by gradient search might not find the actual maximum, \cite{DoyPJ}.   In contrast to the local phenomena, the function $u(D):=\sigma_{\max}(DMD^{-1})$ does not have any local minimal which are not global, so computing $\inf_{D\in\mathcal{D}_{\mathbb{B}} }\sigma_{\max}(DMD^{-1})$ is a reasonable task, \cite{DoyPJ}.  In general, even $\mu_\mathbb{B}(M)\leq \inf_{D\in\mathcal{D}_{\mathbb{B}} }\sigma_{\max}(DMD^{-1})$, it is shown in \cite{DoyPJ} that the equality always holds when $2r+s\leq 3$.   A conditions making the structured singular value reach to the upper bound can be found in Theorem 4.5 in \cite{DoyPJ} and a condition making the structured singular value equal to the spectral norm can also be found in Theorem 4.4 in \cite{DoyPJ}.  

By using these results, Shigeru Yamamoto and Hidenori Kimura showed in \cite{SYHK} that $\mu_\mathbb{B}(M)=\sigma_{\max}(M)$ when the matrix $M$ is a square complex matrix $M$ satisfying $SM=M^TS$ for some signature matrix $S$ and the perturbation block structure $\mathbb{B}$ has $r=0$ or $r=1,s=1$ (with some conditions).  It is well known that if $M$ is a normal matrix, then $\rho(M)=\sigma_{\max}(M)$ and thus by the property (3) mentioned as the earlier, $\rho(M)=\mu_\mathbb{B}(M)=\sigma_{\max}(M)$ for any block structure $\mathbb{B}\subseteq M_n(\mathbb{C})$.  Also, if $M$ is a doubly stochastic matrix, then $\rho(M)=\sigma_{\max}(A)=1$ (this fact does not necessary hold for generalized doubly stochastic matrices) and again $\mu_{\mathbb{B}}(M)=1$ for any block structure $\mathbb{B}\subseteq M_n(\mathbb{C})$.  However, besides these three classes of matrices, according to the best of our knowledge, there is no explicit structured singular values on other classes of matrices.   Most of the literature in the computation area are concerning with the approximation of the upper bound; see, for example, \cite{SiamMutti} and the references therein.

In this article, we investigate and construct a class of generalized row (column) stochastic matrices (square complex matrices having constant row (column) sum) appearing in (\ref{classofexplicitssv}).  Any matrix in this class needs not be normal nor radial (a matrix having the spectral norm equal to the spectral radius).  The explicit formula for the structured singular values of matrices in this class is given in term of the constant row (or column) sum.
 
\section{The results}
Through out this article, we denote the block type structure $\mathbb{B}$ to be a subset of $M_n(\mathbb{C})$ having the form:
$$\mathbb{B}:=\{\operatorname{diag}(\delta_1I_{k_1},\dots,\delta_rI_{k_r},\Delta_1,\dots,\Delta_s)\,|\;\delta_i's\in \mathbb{C}, \Delta_j's\in M_{n_i}(\mathbb{C})\}  $$
for some $r,s\in \mathbb{N}_0$ (set of all nonnegative integers).   Let $A$ be an $n\times n$ generalized row stochastic matrix; namely there exist $c\in \mathbb{C}$ such that $$\sum_{j=1}^n A_{ij}=c,\quad \forall i=1,\dots,n,$$  and denote the largest singular values of $A$ to be $\sigma$; i.e., $\sigma_{\max}(A)=\sigma=\| A \|_2$. We have that:
\begin{prop}\label{mainthm1}
	Let $A\in M_n(\mathbb{C})$ be a generalized row stochastic matrix with the constant row sum $c\in \mathbb{C}$.  Suppose that $\|A\|_2=\sigma$.  Then the following statements hold:
	\begin{enumerate}
		\item $|c|\leq \sigma$.
		\item If  $|c|=\sigma$, then $A$ must have the constant column sum $c\in \mathbb{C}$; namely $A$ must be a $c$-generalized doubly stochastic matrix. 
		\item If $|c|=\sigma$, then $\mu_\mathbb{B}(A^m)=\sigma^m$ for all $m\in \mathbb{N}$ and any block structure $\mathbb{B}\subseteq M_n(\mathbb{C})$.
	\end{enumerate}
\end{prop}
\begin{proof}
	Let $\vec{1}\in \mathbb{C}^n$ be the vector whose all components are one. Since $A\in M_n(\mathbb{C})$ is a generalized row stochastic matrix with the constant row sum $c\in \mathbb{C}$, $$ \sum_{j=1}^n A_{ij}=c,\quad \forall i=1,\dots,n; $$ namely $A\vec{1}=c\vec{1}$. By using the standard 2-norm $\|\cdot\|$ on $\mathbb{C}^n$ and the Cauchy-Schwarz inequality, we have that:
	\begin{equation}\label{CS1}
	 |c|^2\|\vec{1}\|^2=\langle c\vec{1},c\vec{1} \rangle= \langle c\vec{1},A\vec{1} \rangle = \langle cA^* \vec{1},\vec{1} \rangle=|\langle cA^* \vec{1},\vec{1} \rangle |  \leq \|cA^*\vec{1}\|\|\vec{1}\|, 
	 \end{equation}
	where $A^*$ is the conjugate transpose of $A$.  Note that $$\|M\|_2:=\sup_{x\in (\mathbb{C}^*)^n}\frac{\|Mx\|}{\|x\|}\geq \frac{\|My\|}{\|y\|},$$
	for any $M\in M_n(\mathbb{C})$ and any non-zero vector $y\in \mathbb{C}^n$.  Then $\|cA^*\vec{1}\|\leq \|cA^*\|_2\|\vec{1}\|=|c|\|A^*\|_2\|\vec{1}\|$ .  Note also  that $\|A^*\|_2=\|A\|_2$ for any $A\in M_n(\mathbb{C})$.  The inequality (\ref{CS1}) becomes
	\begin{equation}\label{CS2}
|c|^2\|\vec{1}\|^2\leq \langle c A^* \vec{1},\vec{1} \rangle  \leq |c|\|A\|_2\|\vec{1}\|^2=\sigma|c|\|\vec{1}\|^2
	\end{equation}
	By the inequality (\ref{CS2}), if $c\neq 0$, then $\sigma\geq |c|$.  Since $\sigma$ is always non-negative real number, if $c=0$, then $\sigma\geq 0=|c|$.  Hence, the statement (1) holds.
	
	If $0=|c|=\sigma:=\|A\|_2$, then $A=O$ (zero matrix) and thus the statement (2) is obviously true.  Suppose that $0\neq |c|=\sigma$.  Then the inequality (\ref{CS2}) and hence the Cauchy-Schwarz inequality become equalities.  This happens when $cA^*\vec{1}=d'\vec{1}$ for some $d'\in \mathbb{C}$; that is $A^T\vec{1}=d\vec{1}$ where $d:=\bar{d'}/\bar{c}\in \mathbb{C}$; namely,   $$ \sum_{i=1}^n A_{ij}=d,\quad \forall j=1,\dots,n. $$
	By the equality of (\ref{CS2}), we have $$d'\|\vec{1}\|^2=\langle d'\vec{1},1\rangle=\langle c A^* \vec{1},\vec{1} \rangle  = |c|\sigma\|\vec{1}\|^2,  $$
	which means that $d'=|c|\sigma=|c|^2$. This implies that $$d=\frac{\bar{d'}}{\bar{c}}=\frac{|c|^2}{\bar{c}}=c.$$ 
 	Hence, the statement (2) holds.
	
	To prove (3), for the case $|c|=\sigma=\|A\|_2$, we write $c=\sigma e^{i\theta}$ for some $\theta\in \mathbb{R}$.  Then 
	$$ 	A\vec{1}=c\vec{1}=\sigma e^{i\theta} \vec{1};$$
	namely $\sigma e^{i\theta}$ is an eigenvalue of $A$.  Thus,
	$$ \sigma=|\sigma e^{i\theta}|\leq \rho(A)\leq \|A\|_2=\sigma $$
and hence $\rho(A)=\sigma$.  Note that, for each  $m\in \mathbb{N}$, $\rho(A^m)=(\rho(A))^m$ and $\|\cdot\|_2$ is a sub-multiplicative norm; $\|A^m\|_2\leq \|A\|_2^m$.  Then, by using the fact that $\rho(M)\leq \mu_\mathbb{B}(M)\leq \|M\|_2$ for any block structure $\mathbb{B}\subseteq M_n(\mathbb{C})$ and matrix $M\in M_n(\mathbb{C})$, we compute that: $$ \sigma^m=(\rho(A))^m= \rho(A^m)\leq \mu_\mathbb{B}(A^m)\leq \|A^m\|_2\leq\|A\|_2^m=\sigma^m, $$
	for each natural number $m$.  Therefore the statement (3) is proved.
\end{proof}
A consequence of the statement (3) in the above proposition is that:
\begin{cor}
	Let $D_1,D_2,\dots,D_k$ be doubly stochastic matrices of size $n\times n$ and $\delta_1,\delta_2,\dots,\delta_k$ be positive real numbers.  Let $$S_k(D):=\delta_1D_1+\delta_2D_2+\cdots+\delta_kD_k.$$  Then, for each $m\in \mathbb{N}$ and any block structure $\mathbb{B}\subseteq M_n(\mathbb{C})$ $$\mu_\mathbb{B}(S_k(D)^m)=(\delta_1+\delta_2+\cdots+\delta_k)^m.$$ 
\end{cor}
\begin{proof} Note that $S_k(D)$ is a $\delta$-generalized doubly stochastic matrices with constant row sum $\delta:=\delta_1+\delta_2+\cdots+\delta_k$.   Since $S_k(D)\vec{1}=\delta\vec{1}$, $\rho(S_k(D))\geq \delta$.
Note also that $\|D_i\|_2=1$, for each $i=1,\dots,k$.  By triangle inequality for the spectral norm $\|\cdot\|_2$, $$\delta\leq\rho(S_k(D))\leq\|S_k(D)\|_2=\|\delta_1D_1+\delta_2D_2+\cdots+\delta_kD_k\|_2\leq \sum_{i=1}^k\delta_i\|D_i\|_2=\delta.  $$
 Hence, $\|S_k(D)\|_2=\delta$ and thus the result follows from the statement (3) of the above theorem.
\end{proof}
Moreover, in order to calculate the structured singular value, the constant condition on the row sum in the proposition can be relaxed.
\begin{cor}\label{cor2}
	Let $A\in M_n(\mathbb{C})$ be a matrix such that $\|A\|_2=\sigma$.  Let $r_i:=\sum_{k=1}^n A_{ik}$ be the $i$-th row sum for each $i=1,\dots,n$ and $c_j:=\sum_{k=1}^n A_{kj}$ be the $j$-th column sum for each $j=1,\dots,n$.  If $|r_1|=\cdots=|r_n|=\sigma$ or $|c_1|=\cdots=|c_n|=\sigma$, then $$\mu_\mathbb{B}(A^m)=\sigma^m$$ for each $m\in \mathbb{N}$ and any block structure $\mathbb{B}\subseteq M_n(\mathbb{C})$ containing respectively:
	$$\operatorname{diag}(e^{i\theta_1},e^{i\theta_2},\dots,e^{i\theta_n}) \hbox{ or } \operatorname{diag}(e^{i\gamma_1},e^{i\gamma_2},\dots,e^{i\gamma_n}),$$  where $r_j:=\sigma e^{i\theta_j}$ and $c_j:=\sigma e^{i\gamma_j}$ for each $j=1,\dots,n$.
\end{cor}
\begin{proof}
	Note that, for each $M,N\in M_n(\mathbb{C})$, we have that $\rho(M^*)=\rho(M)=\rho(M^T)$ and $\rho(MN)=\rho(NM)$.  Also, by (\ref{lowerexact}), $\mu_\mathbb{B}(M)=\max_{U\in \mathcal{U}_\mathbb{B}}\rho(UM)$.  Then, 
	\begin{equation}\label{Cmu1}
	\mu_\mathbb{B}(M^T)=\mu_\mathbb{B}(M)=\mu_\mathbb{B}(UM)=\mu_\mathbb{B}(MU),  
	\end{equation} 
		for each $U\in \mathcal{U}_\mathbb{B}$. 
	
	 If $|r_1|=\cdots=|r_n|=\sigma$, we then write $r_j=\sigma e^{i\theta_j}$, where $\theta_j\in \mathbb{R}$ for each $j=1,\dots, n$.  Then $A=UA_0$ where $U=\operatorname{diag}(e^{i\theta_1},e^{i\theta_2},\dots,e^{i\theta_n})$ and $A_0$ is a  generalized row stochastic matrix with the constant row sum $\delta\in \mathbb{C}$.  By the condition of the perturbation block structure $\mathbb{B}$, we conclude that $A_0=U^* A$ with $U^*\in \mathcal{U}_\mathbb{B}$.  Hence, by (\ref{Cmu1}), $\mu_\mathbb{B}(A_0)=\mu_\mathbb{B}(A)$ and thus, by the statement (3) in the proposition, we conclude the result.
	 
	 If $|c_1|=\cdots=|c_n|=\sigma$, we then write $c_j=\sigma e^{i\gamma_j}$, where $\gamma_j\in \mathbb{R}$ for each $j=1,\dots, n$.  Then $A^T=WB_0$ where $W=\operatorname{diag}(e^{i\gamma_1},e^{i\gamma_2},\dots,e^{i\gamma_n})$ and $B_0$ is a  generalized row stochastic matrix with the constant row sum $\sigma\in  \mathbb{C}$.  Then $B_0=W^* A^T$ with (by the assumption of $\mathbb{B}$) $W^*\in \mathcal{U}_\mathbb{B}$.  Hence, by (\ref{Cmu1}), $\mu_\mathbb{B}(B_0)=\mu_\mathbb{B}(A)$ and thus by the statement (3) in the proposition, we again conclude the result.
\end{proof}

Note that when $r=0$ or $k_1=k_2=\cdots=k_r=1$, the assumption of the perturbation block structure $\mathbb{B}$ in Corollary \ref{cor2} always be full-filled. Furthermore, any matrix satisfying the assumptions of the above corollary must have the following property.
\begin{prop}\label{exprop}
Let $A\in M_n(\mathbb{C})$ be a matrix such that $\|A\|_2=\sigma$.  Let $r_i:=\sum_{k=1}^n A_{ik}$ be the $i$-th row sum for each $i=1,\dots,n$ and $c_j:=\sum_{k=1}^n A_{kj}$ be the $j$-th column sum for each $j=1,\dots,n$.  The following statements hold:
\begin{enumerate}
	\item  $|r_1|=\cdots=|r_n|=\sigma$ if and only if $A=\sigma W_\theta D$, where $W_\theta:=\operatorname{diag}(e^{i\theta_1},e^{i\theta_2},\dots,e^{i\theta_n})$, $r_k:=\sigma e^{i\theta_k}$, for each $k=1,\dots,n$ and $D$ is a $1$-generalized doubly stochastic matrix with $\|D\|_2=1$.
	\item  $|c_1|=\cdots=|c_n|=\sigma$ if and only if $A=\sigma DW_\gamma $, where $W_\gamma:=\operatorname{diag}(e^{i\gamma_1},e^{i\gamma_2},\dots,e^{i\gamma_n})$, $c_k:=\sigma e^{i\gamma_k}$, for each $k=1,\dots,n$ and $D$ is a $1$-generalized doubly stochastic matrix with $\|D\|_2=1$.
\end{enumerate}
\end{prop}
\begin{proof}
	 Suppose that $|r_1|=\cdots=|r_n|=\sigma$.  We now write $r_k=\sigma e^{i\theta_k}$ for each $k=1,\dots,n$.  Hence $A=\sigma W_\theta D$ such that $D$ is a generalized row stochastic matrix with constant row sum $1$.  Since $\sigma=\|A\|_2=\|\sigma W_\theta D\|_2=\sigma\|D\|_2$, $\|D\|_2=1$.  By the statement (2) in the above proposition, $D$ must be a $1$-generalized doubly stochastic matrix.  The converse is obvious. The statement (2) is proved simply by using the fact that $\|A^T\|_2=\|A\|_2=\sigma$ and by applying $A^T$ to the statement (1).
\end{proof}

By Proposition \ref{exprop}, any matrix in the form $$\sigma \operatorname{diag}(e^{i\theta_1},e^{i\theta_2},\dots,e^{i\theta_n}) D,$$ where $D$ is a doubly stochastic matrix and $\theta_1,\theta_2,\dots,\theta_n\in \mathbb{R}$, $\delta\in \mathbb{R}^+$, will satisfy the assumption of Corollary \ref{cor2}.  For example, the matrix $$A:=\begin{bmatrix}
\frac{i}{20}&0&\frac{i}{20} \\
-\frac{i}{20}&-\frac{i}{20}&0 \\
0&\frac{1}{40}+\frac{i\sqrt{3}}{40}&\frac{1}{40}+\frac{i\sqrt{3}}{40}
\end{bmatrix}  $$ can be written in the form $$ A=\frac{1}{10} \operatorname{diag}(e^{i(\pi /2)},e^{i(-\pi /2)},e^{i(\pi /3)})\begin{bmatrix}
1/2&0&1/2\\
1/2&1/2&0\\
0&1/2&1/2
\end{bmatrix},$$ where the last matrix is a doubly stochastic matrix. By corollary \ref{cor2}, $\mu_{\mathbb{B}}(A^m)=(0.1)^m$ for any $m\in \mathbb{N}$ and any block structure $\mathbb{B}$ containing $\operatorname{diag}(e^{i(\pi /2)},e^{i(-\pi /2)},e^{i(\pi /3)})$.  We also observe that the largest (in absolute value) eigenvalue of $A$ is $0.0488352+0.0447302i$ which means that $\rho(A)<\|A\|_2=0.1$; namely $A$ is not a radial matrix and hence not normal.  

Note further that if $D$ is a $1$-generalized doubly stochastic matrix with $\|D\|_2=1$, then $|D_{ij}|\leq 1$ for each $i,j\in \{1,2,\dots,n\}$. To see this, note that if $\|D\|_2=1$, then $\|Dx\|\leq \|x\|$ for any non-zero vector $x\in \mathbb{C}^n$.  So, if $|D_{ij}|>1$, by applying $x=e_i:=[0,\dots,1,\dots,0]^T\in \mathbb{C}^n$ ($1$ appears only in the $i$-position), we have $$\|Dx\|^2=|D_{1j}|^2+\cdots+|D_{ij}|^2+\cdots+|D_{nj}|^2\geq |D_{ij}|^2>1=\|x\|^2,  $$
which is a contradiction.  Moreover, if $D$ is a $1$-generalized doubly stochastic matrix, then $\|D\|_2$ need not be $1$. For example, the checkerboard matrix
$$D:=\begin{bmatrix}
1&-1&1&-1&1\\
-1&1&-1&1&-1\\
1&-1&1&-1&1\\
-1&1&-1&1&-1\\
1&-1&1&-1&1\\
\end{bmatrix},  $$
is a $1$-generalized doubly stochastic matrix with $\|D\|_2=5$.  Precisely, note first that $Dv=5v$ with $v=[1,-1,1,-1,1]^T$; i.e., $\rho(D)\geq 5$. So, $$ 5\leq \rho(D)\leq \|D\|_2\leq \|D\|_F:=\sqrt{\sum_{i=1}^5\sum_{j=1}^5|D_{ij}|^2}=5. $$  However, there are still many 1-generalized doubly stochastic matrices having spectral norm $1$ which is not a doubly stochastic as we can see some of them from the following theorem.
\begin{thm}
Let $k$ be a positive integer.  Let $a,\alpha_1,\alpha_2,\dots,\alpha_k$ be positive real numbers and $b$ be a real number.
\begin{enumerate}\label{therem2}
	\item If $n=2k$ and $a\geq |b|$, then the circulant matrix $C^e(a,b,\alpha)$ having the row vector $$[a+bi,\; a-bi,\; \alpha_2(a+bi),\; \alpha_2(a-bi),\;\dots,\;\alpha_k(a+bi),\;\alpha_k(a-bi)]$$ as the first row, will be a $\delta_e$-generalized doubly stochastic matrix of size $n\times n$ with the spectral norm $\|C^e(a,b,\alpha)\|_2=\delta_e$, where $$\delta_e:=2a(1+\alpha_2+\cdots+\alpha_k).$$
		\item If $n=2k+1$, then the circulant matrix $C^o(a,b,\alpha)$ having the row vector $$[a+bi,\; a-bi,\; \alpha_2(a+bi),\; \alpha_2(a-bi),\;\dots,\;\alpha_k(a+bi),\;\alpha_k(a-bi),\; \alpha_1]$$ as the first row, will be a $\delta_o$-generalized doubly stochastic matrix of size $n\times n$ with the spectral norm $\|C^o(a,b,\alpha)\|_2=\delta_o$, where $$\delta_o:=2a(1+\alpha_2+\cdots+\alpha_k)+\alpha_1.$$
\end{enumerate}
  
\end{thm}
\begin{proof}
	We first consider the case $n=2k$ and $a\geq |b|$.  It is well known that all eigenvalues of the circulant matrix $C^e(a,b,\alpha)$ are:
	$$ \lambda_j= z+ \bar{z}\omega^j+\alpha_2z\omega^{2j}+ \alpha_2\bar{z}\omega^{3j}+\cdots+\alpha_kz\omega^{2k-2}+\alpha_k\bar{z}\omega^{2k-1};\;\; j=0,1,\dots,n-1,$$
	where $z=a+bi$ and $\omega=e^{2\pi i /n}$.  Then, for each $j\in \{0,1,\dots,n-1\}$,
	\begin{eqnarray*}
	|\lambda_j|&=& |(z+\bar{z}\omega^j)+\alpha_2(z+\bar{z}\omega^j)\omega^{2j}+\alpha_3(z+\bar{z}\omega^j)\omega^{4j}+\cdots+\alpha_k(z+\bar{z}\omega^j)\omega^{2(k-1)j}| \\
	&=&(|z+\bar{z}\omega^j|)(|1+\alpha_2\omega^{2j}+\alpha_3\omega^{4j}+\cdots+\alpha_k\omega^{2(k-1)j}|)\\
	&\leq&(|z+\bar{z}\omega^j|)(1+\alpha_2+\alpha_3+\cdots+\alpha_k).
	\end{eqnarray*}
	The triangle inequality becomes equality when there are non-negative real number $x_t$ such that $$ \alpha_t\omega^{2tj}=x_t\cdot 1=x_t,\quad \hbox{ for each } t=1,\dots,k-1. $$
	Since $\alpha_t>0$, $x_t$ is a positive real number and thus $\omega^{2tj}=1$ for each $t=1,\dots,k-1$.  This implies that $$\omega^{2j}+\omega^{4j} +\cdots+\omega^{2(k-1)j}=k-1 $$
	and that $$\omega^{4j}+\omega^{6j} +\cdots+\omega^{2kj}=\omega^{2j}(k-1) $$
	Subtracting the first equation by the second equation and using $\omega^{2kj}=\omega^{nj}=1^j=1$, we obtain that  $\omega^{2j}-1=(k-1)(1-\omega^{2j})$.  If $\omega^{2j}\neq 1$, then $0\leq k-1=-1$, which is a contradiction.  Then the maximum of $|\lambda_j|$'s can be possibly reached when $\omega^{2j}=1$ or when  $j=0$ or $j=n/2$.  If $j=n/2$, then $\omega^j=-1$ and thus $$|\lambda_{n/2}|=(|z-\bar{z})(1+\alpha_2+\alpha_3+\cdots+\alpha_k)=2|b|(1+\alpha_2+\alpha_3+\cdots+\alpha_k).$$
	If $j=0$, then $\omega^j=1$ and thus $$|\lambda_{0}|=(|z+\bar{z})(1+\alpha_2+\alpha_3+\cdots+\alpha_k)=2a(1+\alpha_2+\alpha_3+\cdots+\alpha_k):=\delta_e.$$
	Since any circulant matrix is always normal and $a\geq |b|$, we conclude that $$\|C^e(a,b,\alpha)\|_2=\max\{|\lambda_j|\;|\;j=0,1,\dots,n-1\}=|\lambda_0|=\delta_e.$$
	
	For the case $n=2k+1$, we compute that: for each $j=0,1,\dots,n-1$, 
		\begin{eqnarray*}
			|\lambda_j|
			&=&|(z+\bar{z}\omega^j)(1+\alpha_2\omega^{2j}+\alpha_3\omega^{4j}+\cdots+\alpha_k\omega^{2(k-1)j})+\alpha_1\omega^{2kj}|\\
			&\leq & |(z+\bar{z}\omega^j)(1+\alpha_2\omega^{2j}+\alpha_3\omega^{4j}+\cdots+\alpha_k\omega^{2(k-1)j})|+|\alpha_1\omega^{2kj}|\\
			&\leq&(|z+\bar{z}\omega^j|)(1+\alpha_2+\alpha_3+\cdots+\alpha_k)+\alpha_1.
		\end{eqnarray*}
By using the same process as the even case, the inequalities become equalities when $j=0$ or $j=n/2$.  Since $n$ is odd, $j=0$.	 Hence,
$$\|C^o(a,b,\alpha)\|_2=\max\{|\lambda_j|\;|\;j=0,1,\dots,n-1\}=|\lambda_0|=2a(1+\alpha_2+\alpha_3+\cdots+\alpha_k)+\alpha_1:=\delta_o.$$	
\end{proof}

  For example, if $n=3$ and $a=1/20, \alpha_1=9/10$ and $b=-\sqrt{3}/20$, then the matrix 
$$C^o:=\begin{bmatrix}
\frac{1}{10}(\frac{1}{2}-\frac{\sqrt{3}i}{2})&\frac{1}{10}(\frac{1}{2}+\frac{\sqrt{3}i}{2})&\frac{9}{10}\\
\frac{9}{10}&\frac{1}{10}(\frac{1}{2}-\frac{\sqrt{3}i}{2})&\frac{1}{10}(\frac{1}{2}+\frac{\sqrt{3}i}{2}) \\ \frac{1}{10}(\frac{1}{2}+\frac{\sqrt{3}i}{2})&
\frac{9}{10} &\frac{1}{10}(\frac{1}{2}-\frac{\sqrt{3}i}{2}) 
\end{bmatrix}, $$
is a $1$-generalized doubly stochastic (which is of course not doubly stochastic matrix) having $1,1,0.7$ as it singular values.  Then, by Corollary \ref{cor2}, $\mu_{\mathbb{B}}((C^o)^m)=1$ for any block structure $\mathbb{B}\subseteq M_3(\mathbb{C})$, and also $$\mu_\mathbb{B}((\sigma W_\theta C^o)^m)=\sigma^m$$ for any $m\in \mathbb{N}$, $\sigma\in \mathbb{R}^+$, $W_\theta:=\operatorname{diag}(e^{i\theta_1},e^{i\theta_2},e^{i\theta_3})\in M_3(\mathbb{C})$ and any block structure $\mathbb{B}$ containing $W_\theta$.  If $n=4$ and $a=1,b=-1/2$ and $\alpha_2=1/3$, we have that $a>|b|$ and thus
$$C^e:=\begin{bmatrix}
1-\frac{i}{2}&1+\frac{i}{2}&\frac{1}{3}-\frac{i}{6}&\frac{1}{3}+\frac{i}{6}\\
\frac{1}{3}+\frac{i}{6}&1-\frac{i}{2}&1+\frac{i}{2}&\frac{1}{3}-\frac{i}{6}\\
\frac{1}{3}-\frac{i}{6}&\frac{1}{3}+\frac{i}{6}&1-\frac{i}{2}&1+\frac{i}{2}\\
1+\frac{i}{2}&\frac{1}{3}-\frac{i}{6}&\frac{1}{3}+\frac{i}{6}&1-\frac{i}{2}\\
\end{bmatrix},  $$
will be an $8/3$-generalized doubly stochastic matrix with $\|C^e\|_2=8/3$.  The condition in statement (1) of Theorem \ref{therem2} is essential.  For example, if we choose $n=4$, $a=1,b=2$  and $\alpha_ 2=3$, then $a<|b|$ and the matrix
$$E:=\begin{bmatrix}
	1+2i&1-2i&3+6i&3-6i\\
3-6i&	1+2i&1-2i&3+6i\\
	3+6i&3-6i&	1+2i&1-2i\\
	1-2i&3+6i&3-6i&	1+2i\\
\end{bmatrix}  $$
is an $8$-generalized doubly stochastic matrix having $\|E\|_2=16$.

Note that any sum or product of circulant matrices is a circulant matrix.  However, a sum of a circulant matrix and a doubly stochastic matrix need not be normal (and hence need not be circulant).  For example, the matrix
$$S:=C^e+\begin{bmatrix}
\frac{1}{4}&0&\frac{3}{8}&\frac{3}{8}\\
	\frac{1}{4}&	0&\frac{3}{8}&\frac{3}{8}\\
\frac{1}{4}&\frac{1}{2}&\frac{5}{32}&\frac{3}{32}\\
\frac{1}{4}&\frac{1}{2}&\frac{3}{32}&\frac{5}{32}\\
\end{bmatrix}= \begin{bmatrix}
\frac{5}{4}-\frac{i}{2}&1+\frac{i}{2}&\frac{17}{24}-\frac{i}{6}&\frac{17}{24}+\frac{i}{6}\\
\frac{7}{12}+\frac{i}{6}&1-\frac{i}{2}&\frac{11}{8}+\frac{i}{2}&\frac{17}{24}-\frac{i}{6}\\
\frac{7}{12}-\frac{i}{6}&\frac{5}{6}+\frac{i}{6}&\frac{37}{32}-\frac{i}{2}&\frac{35}{32}+\frac{i}{2}\\
\frac{5}{4}+\frac{i}{2}&\frac{5}{6}-\frac{i}{6}&\frac{41}{96}+\frac{i}{6}&\frac{37}{32}-\frac{i}{2}\\
\end{bmatrix}, $$
is an $11/3$-generalized doubly stochastic matrix with $\|S\|_2=11/3$, but it is not normal as we can see by a direct calculation that $S^*S\neq SS^*$ .

 Furthermore, it should be remarked that the conclusion of Corollary \ref{CS2} can be extended to the following:
\begin{cor}\label{cor8}
	Let $\tilde{D}_1,\tilde{D}_2,\dots,\tilde{D}_k$ be respectively $\delta_i$-generalized row (or column) stochastic matrices in $M_n(\mathbb{C})$ for which $\|\tilde{D}_i\|_2=\delta_i$ for each positive real number $\delta_1,\dots,\delta_k$.  Let $$S_k(\tilde{D}):=\tilde{D}_1+\tilde{D}_2+\cdots+\tilde{D}_k.$$  Then, $S_k(\tilde{D})$ is a $\delta$-generalized doubly stochastic matrix in $M_n(\mathbb{C})$ such that $\|S_k(\tilde{D})\|_2=\delta$, where $\delta:=\delta_1+\delta_2+\cdots+\delta_k$, and for each $m\in \mathbb{N}$ and any block structure $\mathbb{B}\subseteq M_n(\mathbb{C})$ $$\mu_\mathbb{B}(S_k(\tilde{D})^m)=\delta^m.$$ 
\end{cor}
\begin{proof}
	By Theorem \ref{mainthm1} (statement (2)), $D_i$'s must be $\delta_i$-generalized doubly stochastic matrices. The remaining of the proof is exactly the same arguments as in the proof of Corollary \ref{CS2}. 
\end{proof}
Let $DC(n)$ be  the set of all $n\times n$-doubly circulant matrices in Theorem \ref{therem2} and $D(n)$ be the set of all doubly stochastic matrices.  Let $\Omega_{D}^{Cir}(n)$ be the set of all $n\times n$ complex matrices in the form: $$\Omega_{D}^{Cir}(n):=\{\sum_{i=1}^sd_iD_i+\sum_{j=1}^t\alpha_jC_j\;|\; D_i\in D(n), C_j\in DC(n), d_i,\alpha_i\in \mathbb{R}^+_0, s,t\in \mathbb{N}_0\},$$
where $\mathbb{R}^+_0$ denotes the set of all nonnegative real numbers.  
Then, by Corollary \ref{cor8}, any matrix in $\Omega_{D}^{Cir}(n)$ satisfies the assumptions of Corollary \ref{cor2} and hence the structured singular value, spectral norm and spectral radius are equal to the constant row sum or column sum which are independent from the perturbation block structure.   As the above discussion, matrices in this class are generalized doubly stochastic having spectral norm equal to the constant row sum that are radial but  need not be normal. 

 In conclusion, the structured singular values can be calculated explicitly on the  class of matrices (which need not be radial nor normal):

\begin{equation}\label{classofexplicitssv}
\Omega^m_n(\theta,\gamma,\delta):=\{\delta(W_\theta XW_\gamma)^m\;|\; X\in\Omega_{D}^{Cir}(n); W_\theta,W_\gamma\in Diag(n); m\in \mathbb{N},\delta\in \mathbb{R}^+ \},  
\end{equation}
where $Diag(n):=\{\operatorname{diag}(e^{i\theta_1},e^{i\theta_2},\dots,e^{i\theta_n})\;|\; \theta_1,\dots,\theta_n\in \mathbb{R}\}\subseteq M_n(\mathbb{C})$; namely, by  Corollary \ref{cor2},  if $X\in \Omega_{D}^{Cir}(n)$ has the constant row (column) sum $r$ , then $$\mu_\mathbb{B}(\delta( W_\theta XW_\gamma)^m)=\delta\cdot r^m$$ for any $W_\theta, W_\gamma\in Diag(n) $, $\delta\in \mathbb{R}^+$ , $m\in \mathbb{N}$ and any complex block structure $\mathbb{B}\subseteq M_n(\mathbb{C})$ containing $W_\theta$ and $W_\gamma$.

\section*{Acknowledgments}
The authors would like to thank anonymous referee(s) for reviewing this manuscript.

\end{document}